\documentclass[11pt]{amsart}
\usepackage{amsmath, amsthm, amsopn, amsfonts}

\def\Im{\textrm{Im}}
 
\newtheorem{theoreme}{Theorem}
\newtheorem{proposition}{Proposition}
\newtheorem{lemme}[proposition]{Lemma}

\newtheorem{definition}[proposition]{Definition}
\newtheorem{remarque}[proposition]{Remark}

\numberwithin{equation}{section}
\numberwithin{proposition}{section}

\begin{document}
\title[Global Strichartz Estimates]{Global Strichartz Estimates for 
nontrapping Geometries: About an Article by H. Smith and C. 
Sogge}


\author{N. Burq}
\address{Math\'ematiques, B\^at. 425, Universit\'e Paris Sud Orsay, 91405 Orsay Cedex, France}
\email{nicolas.burq@math.u-psud.fr} \urladdr{www.math.u-psud.fr/\~{}burq}
\thanks{Part of this work was completed during a stay of the author at the Department of mathematics of the University of California, Berkeley, partially funded by a N.A.T.O. fellowship. I thank these two institutions for their support}
\begin{abstract}{The purpose of this note is to present an alternative 
proof of a result by H.~Smith and C.~Sogge showing that in odd 
dimension of space, local (in time) Strichartz estimates and 
exponential decay of the local energy for solutions to wave equations 
imply global Strichartz estimates. Our proof allows to extend the 
result to the case of even dimensions of space}
\end{abstract}
\maketitle
\section{Introduction}
Consider $\Theta\subset B(0, R)\subset  {\mathbb R} ^d$ a smooth ($C^3$) 
obstacle, $\Omega= \Theta ^c$, $g$ a smooth ($C^2$) metric  on $\overline{ \Omega}$ equal to 
$g_{i,j}(x) = \delta_{i,j}$ for $|x| \geq R$,
and the  wave equation on $\Omega$:
\begin{equation}
\label{eq1.1}\left\{\begin{aligned}
(\partial_{ t}^2- \Delta_{g})u(t,x) &= F(t,x)\\
u\mid_{t=0}= f\in \dot{ H}_{D}^\gamma, &\partial _{t}u\mid_{t=0}= g\in 
\dot{ H}_{D}^{\gamma-1}\\
u\mid_{\partial \Omega}&=0
\end{aligned}\right.  
\end{equation}
where $\Delta_{g}$ is the Laplace Beltrami operator on $\Omega$, and 
$\dot{H}^\gamma$, $\gamma\in {\mathbb R}$ are the homogeneous Sobolev 
spaces associated to the square root of the Laplace operator with 
Dirichlet boundary conditions on $\Omega$,  $\sqrt{- 
\Delta_{D}}$. In~\cite{SmSo00} H. Smith and C. Sogge show 
 that in odd space dimensions, if the perturbation $(g, \Theta)$ is 
assumed to be {\em non trapping}, then  local (in time) Strichartz 
estimates  imply global (in time) Strichartz estimates for the 
solutions of the wave equation. The method of proof relies on the 
exponential decay for the local energy of solutions of the 
corresponding wave equation with compactly supported initial data, a localization 
argument on the wave cone and a result by M. Christ and A. 
Kiselev~\cite{ChKi01} allowing to deduce the inhomogeneous estimates 
from the homogeneous ones. The purpose of this note is to generalize 
this result to the case of even dimensions of space. In this case, 
the exponential decay is no more true and one has to replace it by 
another form of decay. In fact the decay we will use (global $L^2$ 
integrability of the local energy norm for {\em any} initial data, 
{\em i.e.} non necessarily compactly supported) is closely related to 
the exponential decay of the local energy used by H. Smith and C. 
Sogge, since the proof of our decay relies on a weaker form of the 
estimates used to prove the exponential decay in scattering 
theory.\par
More precisely, the result by H. Smith and C. Sogge reads as follows:
\begin{definition}[\protect{\cite{SmSo00}}]
\label{de.1}
We say that $1\leq r,s\leq 2\leq p,q\leq + \infty$ and $\gamma$ are 
admissible if the following two mixed norm estimates hold
\begin{itemize}
\item{ \bf Local Strichartz estimates.} For data $f,g,F$ supported in 
$\{ |x| \leq R\}$, and $u$ solution of~\eqref{eq1.1}, one has
\begin{multline}
\label{eq1.2}\|u\|_{L^p_{t}([0,1]; L^q_{x}( \Omega)}+ \|(u, \partial 
_{t}u\|_{L^\infty_{t}([0,1];\dot{H}^\gamma_{D}\times\dot{H}^{\gamma 
-1}_{D}  }\\
\leq C \left( \| f\|_{\dot{H}^\gamma_{D}}+ \|g\|_{\dot{H}^{\gamma 
-1}_{D}}+ \|F\|_{L^r_{t}([0,1]; L^s_{x}( \Omega))}\right)
\end{multline}
\item{\bf Global Minkovski Strichartz Estimates.}  For data $f,g,F$  
and $u$ solution of~\eqref{eq1.1} in the case $\Omega= {\mathbb R}^d$ 
and $g_{i,j}(x)=\delta_{i,j}$, one has 
\begin{equation}
\label{eq1.3}\|u\|_{L^p_{t}({\mathbb R}; L^q_{x}({\mathbb R}^d))}
\leq C \left( \| f\|_{\dot{H}^\gamma_{D}}+ \|g\|_{\dot{H}^{\gamma 
-1}_{D}}+ \|F\|_{L^r_{t}({\mathbb R}; L^s_{x}( {\mathbb R}^d))}\right)
\end{equation}
\end{itemize}
\end{definition}
In~\cite{SmSo00}, H. Smith and C. Sogge show:
\begin{theoreme}[\protect{\cite[Theorem 1.1]{SmSo00}}]
\label{th1} Assume that $d\geq 3$ is odd, $1\leq r,s\leq 2\leq 
p,q\leq + \infty$ and $\gamma$ are admissible, $p>r$ and $\gamma\leq 
\frac{ d-1} 2$. Then for data $f,g,F$ and $u$ solution of~\eqref{eq1.1}, one has
\begin{equation}
\label{eq1.4}\|u\|_{L^p_{t}({\mathbb R}; L^q_{x}( \Omega)}
\leq C \left( \| f\|_{\dot{H}^\gamma_{D}}+ \|g\|_{\dot{H}^{\gamma 
-1}_{D}}+ \|F\|_{L^r_{t}({\mathbb R}; L^s_{x}( \Omega))}\right)
\end{equation}
\end{theoreme}
In this article we are going to prove:
\begin{theoreme}
\label{th2} Assume $d\geq 2$,  $1\leq r,s\leq 2\leq p,q\leq + \infty$ 
and $\gamma$ are admissible, $p>2$ and $\gamma<d/2$. Assume also if $d=2$ that $\Theta\neq \emptyset$. Then for data $f,g,F$ and $u$ 
solution of~\eqref{eq1.1}, then 
\begin{equation}
\label{eq1.4bis}\|u\|_{L^p_{t}({\mathbb R}; L^q_{x}( \Omega))}
\leq C \left( \| f\|_{\dot{H}^\gamma_{D}}+ \|g\|_{\dot{H}^{\gamma 
-1}_{D}}+ \|F\|_{L^r_{t}({\mathbb R}; L^s_{x}( \Omega))}\right)
\end{equation}
\end{theoreme}
\begin{remarque}
\label{re1} The assumption $\Theta\neq \emptyset$ is used to avoid 
low frequency problems. It could be replaced by a ``$0$ is not a resonance'' 
assumption (see section~\ref{se2}).
\end{remarque}
As explained above, the idea in this article is to replace the 
exponential decay of the local energy by another (weaker) form of 
decay. Section~\ref{se2} is hence devoted to the proof of this decay. 
Then in section~\ref{se3} we prove Theorem~\ref{th2}.
\begin{remarque}
At the final stage of this work, I have learned that the results of this article have been obtained independantly (by a slightly different method) by Jason Metcalfe~\cite{Met02} 
\end{remarque}
\section{Global $L^2$-integrability of the local energy}\label{se2} We begin 
first by recalling the definition of the homogeneous Sobolev spaces 
(in the boundary case) $\dot{H}^s_{D}$.
 Remark that we shall restrict ourselves to the case $s< d/2$, hence 
the multiplication by a smooth function $\xi\in C^\infty_{0}({\mathbb 
R}^d)$ is continuous from $ \dot{H}^s( {\mathbb R}^d)$ to $ H^s( 
{\mathbb R}^d)$ and the two norms are equivalent on functions with 
fixed compact support.
Take $R$ large enough so that $\Theta\subset B(0, R)$ and $g_{i,j}(x) 
= \delta_{ i,j}$ for $|x| \geq R$. Fix $\beta \in C^\infty_{0}( 
{\mathbb R}^d)$ equal to $1$ on $\{|x| \leq R\}$, and define for $s\geq 0$ 
$\dot{H}^s_{D}$ to be the set of functions $f\in H^s_{ \text{ loc}}( 
\Omega)$ such that 
\begin{equation}
\|f \|_{\dot{H}^s_{D}}^2= \|\beta f \|_{{H}^s( \Omega)}^2 + \|(1- \beta)f 
\|_{\dot{H}^s({\mathbb {R}}^d)}^2
\end{equation}
is bounded and which satisfy the following compatibility conditions:
\begin{equation}
\label{eqcomp}\Delta_{g}^j f\mid_{\partial \Omega}=0, \forall j; s- 2 
j > \frac 1 2 .
\end{equation}
For $s<0$, define $\dot{H}^s_{D}= (\dot{H}^{-s}_{D})'$
Remark that on  $\dot{H}^s_{D}$, the norms $\|(- \Delta_{g})^ 
{s/2}f\|_{ L^2}$ and $\|f \|_{\dot{H}^s_{D}}$ are equivalent (it is 
clear for $s$ an even integer and follows in general by interpolation 
and duality). Denote by $\dot{{\mathcal {H}}}^ s= \dot{H}^s_{D}\times 
\dot{H}^{s-1}_{D}$ and ${{\mathcal {H}}}^ s= {H}^s_{D}\times 
{H}^{s-1}_{D}$.\par

In this part we are going to prove:
\begin{theoreme}
\label{th3} Consider $\chi\in C^ \infty_{0}( {\mathbb R}^d)$. Assume 
that $d\geq 2$ and if $d=2$ that $\Theta\neq \emptyset$. Then  For 
data $f,g,F$ such that {\em $F$ is supported in $B_{x}( 0, R) \cap 
\Omega$} and $u$ solution of~\eqref{eq1.1}, one has
\begin{equation}
\label{eq1.4ter}\|(\chi u,\chi \partial _{t}u)\|_{L^2_{t}({\mathbb 
R};\dot{H}^\gamma_{D}\times\dot{H}^{\gamma -1}_{D})}
\leq C \left( \| f\|_{\dot{H}^\gamma_{D}}+ \|g\|_{\dot{H}^{\gamma 
-1}_{D}}+ \|F\|_{L^2_{t}({\mathbb R};\dot{H}^{\gamma -1}_{D})}\right)
\end{equation}
\end{theoreme}
\begin{remarque}
\label{re2.1} Remark that in the result above only $F$ (and neither 
$f$ nor $g$) is assumed to be compactly supported.
\end{remarque}
\begin{remarque}
\label{re2.2} As will appear below, the proof of  this result is closely 
related to the proof of the exponential decay of the local energy 
(for compactly supported data). In fact we will use some estimates on 
the resolvent {\em weaker} than the estimates required to prove this 
exponential decay. The result above can be transposed to the 
framework of Schr{\"o}dinger equations. It gives then a (global in 
time) ``smoothing effect''~\cite{BuGeTz02}.
\end{remarque}
\begin{remarque}
\label{re2.3} For $d=2, \gamma\leq 1/2$ and $g_{i,j}(x)= \delta_{i,j}$, for any $u$ 
solution of~\eqref{eq1.1} with $F=0$, we have  (see 
~\cite[Lemma 2.2]{SmSo00})
\begin{equation}
\label{eq1.4quart}\|(\chi u,\chi \partial _{t}u)\|_{L^2_{t}({\mathbb 
R};\dot{H}^\gamma_{D}\times \dot{H}^{\gamma -1}_{D})}
\leq C \left(  \| f\|_{\dot{H}^\gamma_{D}}+\|g\|_{\dot{H}^{\gamma -1}_{D}}\right)
\end{equation}
\end{remarque}
To prove Theorem~\ref{th3}, we first show that by a $TT^*$ argument 
it is enough to prove a non homogeneous estimate. Indeed denote by 
$A=i\left(\begin{smallmatrix} 0 & -\text{Id}\\ - \Delta & 
0\end{smallmatrix}\right)$ the self adjoint operator on ${\mathcal 
{H}}= H^1_{0}( \Omega) \times L^2( \Omega)$ with domain ${\mathcal 
{D}}= \{ (u,v) \in {\mathcal {H}}\cap H^2( \Omega)\times 
H^1_{0}(\Omega)\}$. The solution of~\eqref{eq1.1} with $F=0$ is given 
by 
\begin{equation}
\label{eq2.4}\begin{pmatrix}u\\ \partial _{t}u\end{pmatrix} = e^ {it 
A} \begin{pmatrix} f\\ g \end{pmatrix}
\end{equation}
Remark that the operator $e^{it A}$ is a group of isometries on 
$\dot{ {\mathcal {H}}}^s_{D}$ (equipped with the norm $\|A^ s f 
\|_{L^2}$). \par
Denote by $T= \chi e^{itA}$. The continuity of $T$ from $\dot{ 
{\mathcal {H}}}^s$ to $L^2( {\mathbb R}_{t}; {\mathcal {H}}^s)$ is 
equivalent to the continuity of $T^*$ from $L^2( {\mathbb R}_{t}; 
{\mathcal {H}}^s)$ to  $\dot{ {\mathcal {H}}}^s$. But
\begin{equation}
T^* G= \int _{s\in {\mathbb R}}e^ {-is A}\chi G(s, \cdot) ds
\end{equation}
and the continuity of $T^*$ is in turn equivalent to the continuity 
of $TT^*$ from $L^2( {\mathbb R}_{t}; \dot{\mathcal {H}}^s)$ to $L^2( 
{\mathbb R}_{t}; {\mathcal {H}}^s)$.  Finally 
\begin{equation}
TT^* G (t, \cdot)=\int _{s\in {\mathbb R}}\chi e^ {i(t-s) A}\chi G(s, 
\cdot) ds = \int _{s<t}+ \int_{t<s}
\end{equation}
And it sufficies to prove the continuity of any of the two operators 
above, say for example the continuity of $\chi TT^*_{1}\chi $ with  
\begin{equation}
TT^*_{1} G (t, \cdot)=\int _{s<t}e^ {i(t-s) A} G(s, \cdot) ds
\end{equation}
Remark that this latter result implies in turn the nonhomogeneous 
part of Theorem~\ref{th3}. Suppose that $G$ is supported in $\{s>0\}$ 
(which by density and translation invariance is possible). Then $U= 
TT^*_{1}G$ satisfies the equation 
\begin{equation}
(i \partial _{t}+ A)U= \chi G(t, \cdot), \ U \mid _{t<0}=0
\end{equation}
The Fourier transforms with respect to $t$ of $U$ and $G$  are 
holomorphic in the half plane $\Im z<0$ (due to the support property) 
and  satisfy there 
\begin{equation}
(A-z) \widehat{U}= \chi \widehat{G}
\end{equation}
Taking $z= x-i\varepsilon$ and letting $\varepsilon>0$ tend to $0$,  Theorem~\ref{th3} is a consequence of the fact that the Fourier 
transform is an isometry on $L^2( {\mathbb R}_{t}; {\mathcal {H}})$ 
if ${\mathcal {H}}$ is a Hilbert space and of
\begin{proposition}
\label{le3.1} Under the assumptions of Theorem~\ref{th2}, the 
resolvent $\chi (A-(x-i\varepsilon)) ^ {-1} \chi$ is for $x\in {\mathbb R}$ and $0<|\varepsilon|<1$, uniformly bounded on ${\mathcal 
{H}}^s_{D}$.
\end{proposition}
We are going to deduce Proposition~\ref{le3.1} from the classical result:
\begin{proposition}
\label{prop2}Suppose that the obstacle $\Theta$ is non trapping. Then 
the resolvent of the operator $\Delta_{D}$, $(-\Delta_{D}- 
\lambda)^{-1}$ (which is analytic in ${\mathbb C}\setminus {\mathbb 
R}^+$) satisfies:
\begin{equation}
\label{eq2.2}\begin{gathered}\forall \chi \in C^ \infty_{0}({\mathbb 
R}^2), \exists C>0; \forall \lambda \in {\mathbb R}, 
0<\varepsilon<<1,\\
\| \chi ( -\Delta_{D}- (\lambda\pm i 
\varepsilon))^{-1}\chi\|_{L^2\rightarrow L^2}\leq \frac C { 1 + 
\sqrt{ |\lambda| }}
\end{gathered}
\end{equation}
\end{proposition} 
\begin{remarque}
\label{re4} Proposition~\ref{prop2} was proven for $|\lambda|>>1$ in greater 
generalities by Lax and Phillips~\cite{LaPh89}, Melrose and 
Sj{\"o}strand~\cite{MeSj78, MeSj82}), Vainberg\cite{Va88}, Tang
Zworski~\cite{TaZw00} (see also~\cite{Bu02} for a self contained 
proof which joined with the results in~\cite{Bu97} would relax the 
smoothness assumption required in these papers to a $C^2$ assumption).
The proof for $|\lambda|<<1$ can be found in~\cite[Annexe B.2]{Bu98} 
(see also~\cite{Va68, Va88}). Remark that in~\cite{Bu98} the Poincar{\'e} 
inequality is used to control the local $L^2$-norm by the local 
$L^2$-norm of the gradient of a function (wich is why $\Theta\neq 
\emptyset$ is required). However, the following inequality 
\begin{equation}
\forall \beta< \frac d 2 -1, \forall u\in C^\infty_{0}( {\mathbb 
R}^d), \|\frac u {r (1+ r^2)^{ \beta/2}}\|_{L^2}\leq C \| \frac { 
\nabla u}{(1+ r^2)^{ \beta/2}}\|_{L^2}
\end{equation} allows to handle the argument for $d\geq 4$ (and the 
result is standard in odd dimensions). 
\end{remarque}
Since 
\begin{equation}
(A-z) ^{-1}= \begin{pmatrix}-z(\Delta_{g}+ z^2)^ {-1}& i (\Delta_{g}+ 
z^2)^ {-1}\\
i \Delta_{g}(\Delta_{g}+ z^2)^ {-1} & -z(\Delta_{g}+ z^2)^ {-1}
\end{pmatrix}
\end{equation}
to prove Proposition~\ref{le3.1} we have to estimate uniformly for any $s\in {\mathbb R}$
\begin{align}
&\|\chi (1+|z|)(\Delta_{g}+ z^2)^ {-1}\chi \|_{H^s_{D}\rightarrow H^s_{D} }\label{eq4001}\\
&\|\chi (\Delta_{g}+ z^2)^ {-1}\chi \|_{H^s_{D}\rightarrow H^{s+1}_{D} }\label{eq4002}\\
&\|\chi \Delta_{g}(\Delta_{g}+ z^2)^ {-1}\chi \|_{H^{s}_{D}\rightarrow H^{s-1}_{D} }\label{eq4003}
\end{align}
Remark that using duality and interpolation it is enough to do this for $s\in {\mathbb N}$. \par
For $s=0$, \eqref{eq4001} is Proposition~\ref{prop2} and
if $u=(\Delta_{g}+ z^2)^ {-1}\chi f$, we have 
\begin{equation}
 \begin{aligned}
\int_{\Omega}\chi f \chi \overline{u}&=\int_{ \Omega} (\Delta+ z^2)u \chi \overline{u}\\
&= \int_{ \Omega} -\chi |\nabla u\|^ 2 -\nabla \chi\cdot \nabla u \overline{u}+ z^2 \chi |u| ^2
 \end{aligned}
\end{equation}
which by Cauchy Schwartz implies (using~\eqref{eq4001} for $s=0$ and a different $\chi$) easily~\eqref{eq4002} for $s=0$. By duality we get also~\eqref{eq4002} for $s=-1$. Applying $ \Delta_{g}$ to $ \chi (\Delta_{g}+ z^2)^ {-1} \chi$ and commuting with $\chi$ we obtain by induction~\eqref{eq4002} for $s\in {\mathbb N}$ (hence by duality and interpolation for any $s$). Applying again $\Delta_{g}$ we deduce~\eqref{eq4003} for any $s$. Finally, remark that 
\begin{equation}
\chi \Delta_{g}(\Delta_{g}+ z^2)^ {-1} \chi= \chi \text{Id}\chi +\chi z^2(\Delta_{g}+ z^2)^ {-1} \chi 
\end{equation}from which we estimate uniformly 
\begin{equation}
\|\chi z^2(\Delta_{g}+ z^2)^ {-1} \chi \|_{H^s_{D}\rightarrow H^{s-1}_{D}}.
\end{equation}
Interpolating this latter inequality with~$\text{\eqref{eq4002}}_{s}$ we get~$\text{\eqref{eq4001}}_{s}$.    
\section{From Energy decay to global Strichartz estimates}\label{se3}
In this section we are going to show how the local (in time) 
Strichartz estimates and the global (in time) $L^2$- energy decay 
imply global (in time) Strichartz estimates. We consider firstly the homogeneous case ($F=0$) and we also assume that $\gamma\leq 1/2$ if $d=2$.\par
Consider $\chi \in C^\infty_{0}( {\mathbb R}^d)$ equal to $1$ on 
$B(0, R)$, $u$ solution of~\eqref{eq1.1} with $F=0$ and $v= \chi (x) u$, $w= (1-\chi) (x) u$ solutions of
\begin{equation}
\label{eq3.1}
\left\{ \begin{aligned}
(\partial_{ t}^2- \Delta_{g})v(t,x) &=- [\Delta_{g}, 
\chi] u\\
v\mid_{t=0}&= \chi f\in \dot{ H}^\gamma\\
\partial _{t}v\mid_{t=0}&= \chi g\in \dot{ H}^{\gamma-1}_{D}\\
v\mid_{\partial \Omega}&=0
\end{aligned}\right.  
\end{equation}
\begin{equation}
\label{eq3.1bis}\left\{\begin{aligned}
(\partial_{ t}^2- \Delta_{g})w(t,x) &=(\partial_{ t}^2- 
\Delta_{0})w(t,x) = - [\Delta_{0}, \chi] u\\
w\mid_{t=0}&= (1-\chi) f\in \dot{ H}^\gamma,\\
\partial _{t}w\mid_{t=0}&= (1-\chi) g\in \dot{ H}^{\gamma-1}
\end{aligned}\right.  
\end{equation}
First we show that $w$ satisfies global Strichartz stimates. Since 
$w$ is a solution of the {\em free} wave equation (because it is supported in the set where $g_{i,j}(x)= \delta_{ i,i}$), the contributions 
of $(1-\chi )f$ and $(1-\chi) g$ are handled 
using~\eqref{eq1.3}. Consequently it sufficies to study 
$\widetilde{w}$ the solution of
\begin{equation}
\label{eq3.5}\left\{\begin{aligned}
(\partial_{ t}^2- \Delta_{0})\widetilde{w}(t,x) &=- [\Delta_{0}, 
\chi] u\\
\widetilde{w}\mid_{t=0}&= 0\\
\partial _{t}\widetilde{w}\mid_{t=0}&=0
\end{aligned}\right.  
\end{equation}
  For this we use the following result by M. Christ and A. 
Kiselev~\cite{ChKi01}:
\begin{theoreme}[M. Christ and A. Kiselev]
\label{thCk} Consider a bounded operator $T:L^p( {\mathbb R}; 
B_{1})\rightarrow L^q( {\mathbb R}; B_{2})$ given by a locally 
integrable kernel $K(t,s)$ with value operators from $B_{1}$ to 
$B_{2}$ where $B_{1;2}$ are Banach spaces. Suppose that $p<q$. Then 
the operator 
\begin{equation}
\widetilde {T} f(t) = \int_{s<t} K(t,s) f(s) ds
\end{equation}
is bounded from $L^p( {\mathbb R}; B_{1})$ to $L^q( {\mathbb R}; 
B_{2})$ by
\begin{equation}
\|\widetilde {T}\|_{L^p( {\mathbb R}; B_{1})\rightarrow L^q( {\mathbb 
R}; B_{2})}\leq   (1-2^{ - p^{-1}- q^{-1}})^{-1} \| T\|_{L^p( {\mathbb R}; 
B_{1})\rightarrow L^q( {\mathbb R}; B_{2})}
\end{equation}
\end{theoreme}
According to~\eqref{eq3.5},
\begin{equation}
(\widetilde{w}, \partial _{t}\widetilde{w})(t,x)= \int_{0}^t 
e^{i(t-s)A_{0}}- [\Delta_{g}, \chi] u(s, \cdot) ds
\end{equation}
According to Theorem~\ref{th3} (and Remark~\ref{re2.3} since $\gamma\leq 1/2$ if $d=2$) applied to 
the free wave operator, for any $\chi \in C^\infty_{0}( {\mathbb 
R}^d)$ the operator $T_{0}= \chi e^{itA_{0}}$ maps $\dot{H}^ 
\gamma\times \dot {H}^ { \gamma-1}$ into $L^2( {\mathbb 
R}_{t}^+;\dot{H}^ \gamma\times \dot {H}^ { \gamma-1})$ and 
consequently its adjoint 
\begin{equation}
T_{0}^* G= \int _{0}^\infty e^{-isA_{0}} \chi G ds
\end{equation}
maps  $L^2( {\mathbb R}_{t}^+;\dot{H}^ \gamma\times \dot {H}^ { 
\gamma-1})$ into  $\dot{H}^ \gamma\times \dot {H}^ { \gamma-1}$ and since
\begin{equation}
\left( \widetilde {w}, \partial _{t}\widetilde {w}\right) = - \int_{0}^t e^{i(t-s)A_{0}}\left( 0, [\Delta_{g}, \chi] u\right)ds,
\end{equation}
the global Strichartz estimate follows from Theorem~\ref{thCk} (here 
we are speaking about adjoints relative to the $L^2$ in time and 
$\dot{H}^ \gamma\times \dot {H}^ { \gamma-1}$ in space norms).
\par We come back to the analysis of $v= \chi u$. Consider 
$\varphi\in C^\infty_{0}(]0, 1[)$ equal to $1$ on $[1/4,3/4]$. Denote 
by $v_{n}=\varphi (t-n/2) v$ solution of
\begin{equation}
\label{eq3.5bis}
\left\{\begin{aligned}
(\partial_{ t}^2- \Delta_{g})v_{n}(t,x) &=-\varphi 
(t-n/2)[\Delta_{g}, \chi] u+ [\partial _{t}^2, \varphi(t-n/2)] \chi u\\
&= u_{n}(t, x)\\
v_{n}\mid_{t<n/2}&= 0\\
\partial _{t}v_{n}\mid_{t<n/2}&=0\\
v_{n\mid_{\partial \Omega}}&=0
\end{aligned}\right.  
\end{equation}
With, according to Theorem~\ref{th3},
\begin{equation}
\label{eq4.2}
\sum_{n\in {\mathbb Z}}\|u_{n}\|^2_{L^2( {\mathbb R}_{t}; \dot{H}^{ 
\gamma-1}_{ D})}\leq C(\|f\|^2_{\dot{H}^\gamma _{D}}+ 
\|g\|^2_{\dot{H}^{\gamma-1} _{D}}) 
\end{equation}
 According to the (local in time) Strichartz inequality~\eqref{eq1.2} 
(and Minkovski inequality), we obtain
\begin{equation}
\label{eq1.2bis}\|v_{n}\|_{L^p_{t}({\mathbb R}_{t}; L^q_{x}( \Omega))}+ \|(v_{n}, 
\partial 
_{t}v_{n}\|_{L^\infty_{t}([0,1];\dot{H}^\gamma_{D}\times\dot{H}^{\gamma 
-1}_{D}  )}
\leq C \left(  \|u_{n}\|_{L^1_{t}([-1+n/2,2+n/2];H^{ \gamma-1}_{x}( 
\Omega))}\right)
\end{equation}
Since 
\begin{equation}
\|v\|^p_{L^p_{t}({\mathbb R}; L^q_{x}( \Omega))}\sim \sum_{n\in 
{\mathbb Z}}\|v_{n}\|^p_{L^p_{t}({\mathbb R}_{t}; L^q_{x}( \Omega)}
\end{equation}
and since  $p\geq 2$, we obtain, using~\eqref{eq4.2} and~\eqref{eq1.2bis}
\begin{equation}
 \begin{aligned}
\|v\|^p_{L^p_{t}({\mathbb R}; L^q_{x}( \Omega))}&\leq \sum_{n\in 
{\mathbb Z}}\|v_{n}\|^p_{L^p_{t}({\mathbb R}_{t}; L^q_{x}( \Omega))} \leq C \sum_{n\in 
{\mathbb Z}}\|u_{n}\|^p_{L^1_{t} ;H^{ \gamma-1}_{x}( 
\Omega)}\\
&\leq C \bigl( \sum_{n\in 
{\mathbb Z}} \|u_{n}\|^2_{L^2_{t};H^{ \gamma-1}_{x}( 
\Omega)}\bigr)^{p/2} \leq C (\|f\|^2_{\dot{H}^\gamma _{D}}+ 
\|g\|^2_{\dot{H}^{\gamma-1} _{D}})^{p/2}
 \end{aligned}
\end{equation}
which proves Theorem~\ref{th2} for $F=0$. As recognized by H. Smith and C. Sogge~\cite{SmSo00}, another application of Theorem~\ref{thCk} gives the general case.\par 
Finally, if $d=2$ we have to get rid of the assumption $\gamma\leq 1/2$. Interpolating between~\eqref{eq1.4quart} and the energy estimate (for any $\gamma<1$)
\begin{equation}
\|(\chi u, \chi \partial _{u})\|_{L^ \infty ({\mathbb R} _{t};\dot{H}^{\gamma} _{D}\times  \dot{H}^{\gamma-1})}\leq \|(f, g)\|_{L^ \infty ({\mathbb R} _{t};\dot{H}^{\gamma} _{D}\times  \dot{H}^{\gamma-1})}
\end{equation}
we obtain for any $q>1/(1-\gamma)$
\begin{equation}
\|(\chi u, \chi \partial _t {u})\|_{L^ q ({\mathbb R} _{t};\dot{H}^{\gamma} _{D}\times  \dot{H}^{\gamma-1})}\leq \|(f, g)\|_{L^ \infty ({\mathbb R} _{t};\dot{H}^{\gamma} _{D}\times  \dot{H}^{\gamma-1})}
\end{equation}
using this to estimate $\widetilde w$ and Christ-Kiselev Theorem we can conclude as above provided $p>1/(1-\gamma)$. Finally we remark that this latter condition is automatically fulfilled due to the scaling condition $\frac 1 p + \frac d q = \frac d 2 - \gamma$ which is necessary for~\eqref{eq1.2} to hold. 

\end{document}